\documentclass{amsart}
\usepackage{amssymb}
\usepackage{amsfonts}
\usepackage{graphicx}
\usepackage{amscd}

\newtheorem{theorem}{Theorem}
\theoremstyle{plain}

\newtheorem{corollary}{Corollary}

\newtheorem{lemma}{Lemma}

\newtheorem{proposition}{Proposition}
\newtheorem{remark}{Remark}

\numberwithin{equation}{section}

\input{tcilatex}

\begin{document}
\title[Bochner's Integral]{Norm Estimates for the Difference Between
Bochner's Integral and the Convex Combination of Function's Values}
\author{P. Cerone}
\address{School of Computer Science and Mathematics\\
Victoria University of Technology\\
PO Box 14428, MCMC 8001\\
VICTORIA, Australia.}
\email{pietro.cerone@vu.edu.au}
\urladdr{http://rgmia.vu.edu.au/cerone/index.html}
\author{Y.J. Cho$^{\bigstar }$}
\address{Department of Mathematics Education\\
College of Education, Gyeongsang National University\\
Chinju 660-701, KOREA}
\email{yjcho@nongae.gsnu.ac.kr}
\author{S.S. Dragomir}
\address{School of Computer Science and Mathematics\\
Victoria University of Technology\\
PO Box 14428, MCMC 8001\\
VICTORIA, Australia.}
\email{sever.dragomir@vu.edu.au}
\urladdr{http://rgmia.vu.edu.au/SSDragomirWeb.html}
\author{J.K. Kim}
\address{Department of Mathematics\\
Kyungnam University, Masan,\\
Kyungnam 631-701, Korea}
\email{jongkyuk@kyungnam.ac.kr}
\author{S.S. Kim$^{\blacklozenge }$}
\address{Department of Mathematics\\
Dongeui University, Pusan 614-714, Korea}
\email{sskim@dongeui.ac.kr}
\date{July 16, 2003.}
\keywords{Bochner's Integral, Ostrowski Inequality, Quadrature Formulae\\
$\ \bigstar ,\blacklozenge $ \ Corresponding authors.}

\begin{abstract}
Norm estimates are developed between the Bochner integral of a vector-valued
function in Banach spaces having the Radon-Nikodym property and the convex
combination of function values taken on a division of the interval $\left[
a,b\right] .$
\end{abstract}

\maketitle

\section{Introduction}

A Banach space $X$ with the property that every absolutely continuous $X-$%
valued function is almost everywhere differentiable is said to be a \textit{%
Radon-Nikodym }space \cite[pp. 217--219]{4ab} or \cite{3b,2b} (see also \cite%
{BBCD}). For example, every reflexive Banach space (in particular, every
Hilbert space) is a Radon-Nikodym space, but the space $L_{\infty }\left[ 0,1%
\right] $ of all $\mathbb{K}-$valued, essentially bounded functions defined
on the interval $\left[ 0,1\right] $, endowed with the norm 
\begin{equation*}
\left\Vert g\right\Vert _{\infty }:=ess\sup\limits_{t\in \left[ 0,1\right]
}\left\vert g\left( t\right) \right\vert ,
\end{equation*}%
is a Banach space which is not a Radon-Nikodym space.

A function $f:\left[ a,b\right] \rightarrow X$ is \textit{measurable }if
there exists a sequence of simple functions $\left( f_{n}\right) $ (with $%
f_{n}:\left[ a,b\right] \rightarrow X$) which converges punctually a.e. on $%
\left[ a,b\right] $.

It is well-known that a measurable function $f:\left[ a,b\right] \rightarrow
X$ is Bochner integrable if and only if its norm, that is, the function $%
t\longmapsto \left\Vert f\right\Vert \left( t\right) :=\left\Vert f\left(
t\right) \right\Vert :\left[ a,b\right] \rightarrow X$ is Lebesgue
integrable on $\left[ a,b\right] $, (see for example \cite{HP}). The Bochner
integral of $f$ shall be represented by $\left( B\right) \int_{a}^{b}f.$

Further, we use the integration by parts formula, which holds under the
following general conditions:

Let $-\infty <a<b<\infty $ and $f,g$ \ be two mappings defined on $\left[ a,b%
\right] $ such that $f$ is $\mathbb{C}$-valued and $g$ is $X$-valued, where $%
X$ is a real or complex Banach space. If $f,g$ are differentiable on $\left[
a,b\right] $ and their derivatives are Bochner integrable on $\left[ a,b%
\right] $, then 
\begin{equation*}
\left( B\right) \int_{a}^{b}f^{\prime }g=f\left( b\right) g\left( b\right)
-f\left( a\right) g\left( b\right) -\left( B\right) \int_{a}^{b}fg^{\prime }.
\end{equation*}

For some results on the Ostrowski inequality for real-valued functions, see 
\cite{0b}, \cite{1a0b}, \cite{1cb} and \cite{1db}, and the references
therein.

The following theorem concerning a version of Ostrowski's inequality for
vector-valued functions has been obtained in \cite{BBCD}.

\begin{theorem}
\label{t2.1}Let $\left( X;\left\| \cdot \right\| \right) $ be a Banach space
with the Radon-Nikodym property and $f:\left[ a,b\right] \rightarrow X$ an
absolutely continuous function on $\left[ a,b\right] $ with the property
that $f^{\prime }\in L_{\infty }\left( \left[ a,b\right] ;X\right) $, i.e., 
\begin{equation*}
\left| \left\| f^{\prime }\right\| \right| _{\left[ a,b\right] ,\infty
}:=ess\sup\limits_{t\in \left[ a,b\right] }\left\| f^{\prime }\left(
t\right) \right\| <\infty .
\end{equation*}
Then we have the inequalities: 
\begin{eqnarray}
&&\left\| f\left( s\right) -\frac{1}{b-a}\left( B\right) \int_{a}^{b}f\left(
t\right) dt\right\| \\
&\leq &\frac{1}{b-a}\left[ \int_{a}^{s}\left( t-a\right) \left\| f^{\prime
}\left( t\right) \right\| dt+\int_{s}^{b}\left( b-t\right) \left\| f^{\prime
}\left( t\right) \right\| dt\right]  \notag \\
&\leq &\frac{1}{2\left( b-a\right) }\left[ \left( s-a\right) ^{2}\left|
\left\| f^{\prime }\right\| \right| _{\left[ a,s\right] ,\infty }+\left(
b-s\right) ^{2}\left| \left\| f^{\prime }\right\| \right| _{\left[ s,b\right]
,\infty }\right]  \notag \\
&\leq &\left[ \frac{1}{4}+\left( \frac{s-\frac{a+b}{2}}{b-a}\right) ^{2}%
\right] \left( b-a\right) \left| \left\| f^{\prime }\right\| \right| _{\left[
a,b\right] ,\infty }  \notag \\
&\leq &\frac{1}{2}\left( b-a\right) \left| \left\| f^{\prime }\right\|
\right| _{\left[ a,b\right] ,\infty };  \notag
\end{eqnarray}
for any $s\in \left[ a,b\right] $, where $\left( B\right)
\int_{a}^{b}f\left( t\right) dt$ is the Bochner integral of $f$.
\end{theorem}

Bounds involving the $p-$norms, $p\in \lbrack 1,\infty ),$ of the derivative 
$f^{\prime }$, are embodied in the following theorem \cite{BBCD}.

\begin{theorem}
\label{ta}Let $\left( X,\left\Vert \cdot \right\Vert \right) $ be a Banach
space with the Radon-Nikodym property and $f:\left[ a,b\right] \rightarrow X$
be an absolutely continuous function on $\left[ a,b\right] $ with the
property that $f^{\prime }\in L_{p}\left( \left[ a,b\right] ;X\right) $, $%
p\in \lbrack 1,\infty )$, i.e., 
\begin{equation}
\left\vert \left\Vert f^{\prime }\right\Vert \right\vert _{\left[ a,b\right]
,p}:=\left( \int_{a}^{b}\left\Vert f^{\prime }\left( t\right) \right\Vert
^{p}dt\right) ^{\frac{1}{p}}<\infty .  \label{a1}
\end{equation}%
Then we have the inequalities: 
\begin{eqnarray}
&&\left\Vert f\left( s\right) -\frac{1}{b-a}\left( B\right)
\int_{a}^{b}f\left( t\right) dt\right\Vert  \label{a2} \\
&\leq &\frac{1}{b-a}\left[ \int_{a}^{s}\left( t-a\right) \left\Vert
f^{\prime }\left( t\right) \right\Vert dt+\int_{s}^{b}\left( b-t\right)
\left\Vert f^{\prime }\left( t\right) \right\Vert dt\right]  \notag
\end{eqnarray}%
\begin{eqnarray*}
&\leq &\left\{ 
\begin{array}{l}
\dfrac{1}{b-a}\left[ \left( s-a\right) \left\vert \left\Vert f^{\prime
}\right\Vert \right\vert _{\left[ a,s\right] ,1}+\left( b-s\right)
\left\vert \left\Vert f^{\prime }\right\Vert \right\vert _{\left[ s,b\right]
,1}\right] \\ 
\;\;\;\;\;\;\;\;\;\;\;\;\;\;\;\;\;\;\;\;\;\;\;\;\;\;\;\;\;\;\;\;\;\;\;\;\;\;%
\;\;\;\;\;\;\;\;\;\;\;\;\;\;\;\;\;\;\;\;\;\;\;\;\;\;\;\text{if \ }f^{\prime
}\in L_{1}\left( \left[ a,b\right] ;X\right) ; \\ 
\\ 
\dfrac{1}{\left( b-a\right) \left( q+1\right) ^{\frac{1}{q}}}\left[ \left(
s-a\right) ^{\frac{1}{q}+1}\left\vert \left\Vert f^{\prime }\right\Vert
\right\vert _{\left[ a,s\right] ,p}+\left( b-s\right) ^{\frac{1}{q}%
+1}\left\vert \left\Vert f^{\prime }\right\Vert \right\vert _{\left[ s,b%
\right] ,p}\right] \\ 
\;\;\;\;\;\;\;\;\;\;\;\;\;\;\;\;\;\;\;\;\;\;\;\;\;\;\;\;\;\;\;\;\;\;\text{if
\ }p>1,\;\frac{1}{p}+\frac{1}{q}=1\;\text{and }f^{\prime }\in L_{p}\left( %
\left[ a,b\right] ;X\right) ;%
\end{array}%
\right. \\
&& \\
&\leq &\left\{ 
\begin{array}{l}
\left[ \dfrac{1}{2}+\left\vert \dfrac{s-\frac{a+b}{2}}{b-a}\right\vert %
\right] \left\vert \left\Vert f^{\prime }\right\Vert \right\vert _{\left[ a,b%
\right] ,1}\;\;\;\;\;\;\;\;\;\;\;\;\;\text{if \ }f^{\prime }\in L_{1}\left( %
\left[ a,b\right] ;X\right) ; \\ 
\\ 
\dfrac{1}{\left( q+1\right) ^{\frac{1}{q}}}\left[ \left( \dfrac{s-a}{b-a}%
\right) ^{q+1}+\left( \dfrac{b-s}{b-a}\right) ^{q+1}\right] ^{\frac{1}{q}%
}\left( b-a\right) ^{\frac{1}{q}}\left\vert \left\Vert f^{\prime
}\right\Vert \right\vert _{\left[ a,b\right] ,p} \\ 
\;\;\;\;\;\;\;\;\;\;\;\;\;\;\;\;\;\;\;\;\;\;\;\;\;\;\;\;\;\;\;\;\;\;\;\;\;\;%
\;\;\;\;\;\;\;\;\;\;\;\;\;\;\;\;\;\;\;\;\;\text{if\ }f^{\prime }\in
L_{p}\left( \left[ a,b\right] ;X\right) .%
\end{array}%
\right.
\end{eqnarray*}
\end{theorem}

The main aim of this paper is to point out estimates between the Bochner
integral of a vector-valued function, with values in Banach spaces having
the Radon-Nikodym property and a convex combination of values taken on a
given division of the interval $\left[ a,b\right] .$ The obtained results
naturally extend the Ostrowski type inequalities mentioned above. Some
particular cases for two and three points rules are also given.

\section{The Results}

Let $a\leq b$ and $c\in \mathbb{R}$. Define the mapping 
\begin{equation}
\mu _{p}\left( a,c,b\right) :=\left\{ 
\begin{array}{c}
\int_{a}^{b}\left\vert t-c\right\vert ^{p}dt\;\;\;\text{if\ \ \ }p\in
\lbrack 1,\infty ); \\ 
\\ 
\max\limits_{t\in \left[ a,b\right] }\left\vert t-c\right\vert \;\;\;\text{%
if\ \ \ }p=\infty .%
\end{array}%
\right.   \label{2.0}
\end{equation}%
We observe that:

\begin{enumerate}
\item If $c<a,$ then 
\begin{align*}
\mu _{p}\left( a,c,b\right) & =\int_{a}^{b}\left( t-c\right) ^{p}dt \\
& =\frac{1}{p+1}\left[ \left( b-c\right) ^{p+1}-\left( a-c\right) ^{p+1}%
\right] ,\;\;\;\text{for\ \ \ }p\in \lbrack 1,\infty )
\end{align*}%
and 
\begin{equation*}
\mu _{\infty }\left( a,c,b\right) =b-c.
\end{equation*}

\item If $c\in \left[ a,b\right] ,$ then 
\begin{align*}
\mu _{p}\left( a,c,b\right) & =\int_{a}^{c}\left( c-t\right)
^{p}dt+\int_{c}^{b}\left( t-c\right) ^{p}dt \\
& =\frac{1}{p+1}\left[ \left( c-a\right) ^{p+1}+\left( b-c\right) ^{p+1}%
\right]
\end{align*}%
for $p\in \lbrack 1,\infty )$ and 
\begin{equation*}
\mu _{\infty }\left( a,c,b\right) =\max \left( c-a,b-c\right) =\frac{1}{2}%
\left( b-a\right) +\left\vert c-\frac{a+b}{2}\right\vert .
\end{equation*}

\item If $b<c,$ then 
\begin{align*}
\mu _{p}\left( a,c,b\right) & =\int_{a}^{b}\left( c-t\right) ^{p}dt \\
& =\frac{1}{p+1}\left[ \left( c-a\right) ^{p+1}-\left( c-b\right) ^{p+1}%
\right] ,\;\;\;\text{for\ \ \ }p\in \lbrack 1,\infty )
\end{align*}%
and 
\begin{equation*}
\mu _{\infty }\left( a,c,b\right) =c-a.
\end{equation*}
\end{enumerate}

Consequently, we may conclude that 
\begin{equation*}
\mu _{p}\left( a,c,b\right) =\left\{ 
\begin{array}{ll}
\frac{1}{p+1}\left[ \left( b-c\right) ^{p+1}-\left( a-c\right) ^{p+1}\right]
& \text{if \hspace{0.05in}}c<a; \\ 
&  \\ 
\frac{1}{p+1}\left[ \left( c-a\right) ^{p+1}+\left( b-c\right) ^{p+1}\right]
& \text{if \hspace{0.05in}}c\in \left[ a,b\right] ; \\ 
&  \\ 
\frac{1}{p+1}\left[ \left( c-a\right) ^{p+1}-\left( c-b\right) ^{p+1}\right]
& \text{if \hspace{0.05in}}b<c;%
\end{array}%
\right.
\end{equation*}%
for $p\in \lbrack 1,\infty )$ and 
\begin{equation*}
\mu _{\infty }\left( a,c,b\right) =\left\{ 
\begin{array}{ll}
b-c & \text{if \hspace{0.05in}}c<a; \\ 
&  \\ 
\frac{1}{2}\left( b-a\right) +\left\vert c-\frac{a+b}{2}\right\vert & \text{%
if \hspace{0.05in}}c\in \left[ a,b\right] ; \\ 
&  \\ 
c-a & \text{if \hspace{0.05in}}b<c;%
\end{array}%
\right.
\end{equation*}%
where $\mu _{s}\left( a,c,b\right) ,s\in \left[ 1,\infty \right] $ is as
defined in (\ref{2.0}).

The following integral identity is of interest.

\begin{lemma}
\label{l1}Let $f:\left[ a,b\right] \rightarrow X$ be an absolutely
continuous function on the Banach space $X,$ $X$ is with the property of
Radon-Nikodym, $a\leq x_{1}\leq \cdots \leq x_{n-1}\leq x_{n}\leq b$ and $%
p_{i}>0$ $\left( i=1,\dots ,n\right) $ with $\sum_{i=1}^{n}p_{i}=1.$ Then we
have the identity: 
\begin{multline}
\sum_{i=1}^{n}p_{i}f\left( x_{i}\right) -\frac{1}{b-a}\left( B\right)
\int_{a}^{b}f\left( t\right) dt=\frac{1}{b-a}\left( B\right)
\int_{a}^{x_{1}}\left( t-a\right) f^{\prime }\left( t\right) dt  \label{2.1}
\\
+\frac{1}{b-a}\sum_{i=1}^{n-1}\left( B\right) \int_{x_{i}}^{x_{i+1}}\left[
t-\left( P_{i}b+\bar{P}_{i}a\right) \right] f^{\prime }\left( t\right) dt \\
+\frac{1}{b-a}\left( B\right) \int_{x_{n}}^{b}\left( t-b\right) f^{\prime
}\left( t\right) dt,
\end{multline}%
where $\left( B\right) \int_{a}^{b}f\left( t\right) dt$ is the Bochner
integral, $P_{i}:=\sum_{k=1}^{i}p_{k}$ and $\bar{P}_{i}=1-P_{i}.$

The sum in the middle is assumed to be zero when $n=1.$
\end{lemma}

\begin{proof}
We know that, on utilizing the integration by parts formula, for any $x\in %
\left[ a,b\right] ,$ we have the representation (see for example \cite{BBCD}%
) 
\begin{equation}
f\left( x\right) =\frac{1}{b-a}\left( B\right) \int_{a}^{b}f\left( t\right)
dt+\frac{1}{b-a}\left( B\right) \int_{a}^{b}k\left( x,t\right) f^{\prime
}\left( t\right) dt,  \label{2.2}
\end{equation}%
where 
\begin{equation*}
k\left( x,t\right) =\left\{ 
\begin{array}{ll}
t-a & \text{if \hspace{0.05in}}a\leq t\leq x\leq b, \\ 
&  \\ 
t-b & \text{if \hspace{0.05in}}a\leq x<t\leq b.%
\end{array}%
\right.
\end{equation*}%
Putting in (\ref{2.2}) $x=x_{i}$ $\left( i=1,\dots ,n\right) ,$ multiplying
by $p_{i}\geq 0$ and summing over $i$ from $1$ to $n,$ we deduce 
\begin{equation}
\sum_{i=1}^{n}p_{i}f\left( x_{i}\right) =\frac{1}{b-a}\left( B\right)
\int_{a}^{b}f\left( t\right) dt+\frac{1}{b-a}\left( B\right) \int_{a}^{b}%
\left[ \sum_{i=1}^{n}p_{i}k\left( x_{i},t\right) \right] f^{\prime }\left(
t\right) dt.  \label{2.3}
\end{equation}%
However, 
\begin{gather*}
k\left( x_{1},t\right) =\left\{ 
\begin{array}{ll}
t-a & \text{if \hspace{0.05in}}a\leq t\leq x_{1}\leq b, \\ 
&  \\ 
t-b & \text{if \hspace{0.05in}}a\leq x_{1}<t\leq b,%
\end{array}%
\right. \\
\cdots \cdots \cdots \cdots \cdots \cdots \cdots \cdots \cdots \\
k\left( x_{n},t\right) =\left\{ 
\begin{array}{ll}
t-a & \text{if \hspace{0.05in}}a\leq t\leq x_{n}\leq b, \\ 
&  \\ 
t-b & \text{if \hspace{0.05in}}a\leq x_{n}<t\leq b,%
\end{array}%
\right.
\end{gather*}%
then 
\begin{equation}
S\left( \bar{x},\bar{p},t\right) :=\sum_{i=1}^{n}p_{i}k\left( x_{i},t\right)
\label{2.4}
\end{equation}%
\begin{equation*}
=\left\{ 
\begin{array}{ll}
p_{1}\left( t-a\right) +p_{2}\left( t-a\right) +\cdots +p_{n-1}\left(
t-a\right) +p_{n}\left( t-a\right) , & a\leq t\leq x_{1}\leq b, \\ 
p_{1}\left( t-b\right) +p_{2}\left( t-a\right) +\cdots +p_{n-1}\left(
t-a\right) +p_{n}\left( t-a\right) , & a\leq x_{1}<t\leq x_{2}\leq b, \\ 
&  \\ 
\;\;\;\;\;\;\;\;\;\;\;\;\;\;\;\;\;\;\;\;\;\;\;\;\;\;\;\;\;\cdots \cdots
\cdots \cdots \cdots \cdots \cdots \cdots \cdots &  \\ 
&  \\ 
p_{1}\left( t-b\right) +p_{2}\left( t-b\right) +\cdots +p_{n-1}\left(
t-b\right) +p_{n}\left( t-a\right) , & a\leq x_{n-1}\leq t\leq x_{n}\leq b,
\\ 
p_{1}\left( t-b\right) +p_{2}\left( t-b\right) +\cdots +p_{n-1}\left(
t-b\right) +p_{n}\left( t-b\right) , & a\leq x_{n}<t\leq b,%
\end{array}%
\right.
\end{equation*}%
\begin{equation*}
=\left\{ 
\begin{array}{ll}
t-a, & a\leq t\leq x_{1}\leq b, \\ 
p_{1}\left( t-b\right) +\left( p_{2}+\cdots +p_{n}\right) \left( t-a\right) ,
& a\leq x_{1}<t\leq x_{2}\leq b, \\ 
&  \\ 
\;\;\;\;\;\;\;\;\;\;\;\;\;\;\;\;\;\;\;\cdots \cdots \cdots \cdots \cdots
\cdots \cdots \cdots \cdots &  \\ 
&  \\ 
\left( p_{1}+\cdots +p_{n-1}\right) \left( t-b\right) +p_{n}\left(
t-a\right) , & a\leq x_{n-1}\leq t\leq x_{n}\leq b, \\ 
t-b, & a\leq x_{n}<t\leq b,%
\end{array}%
\right.
\end{equation*}%
\begin{equation*}
=\left\{ 
\begin{array}{ll}
t-a, & a\leq t\leq x_{1}\leq b, \\ 
t-\left[ p_{1}b+\left( p_{2}+\cdots +p_{n}\right) a\right] , & a\leq
x_{1}<t\leq x_{2}\leq b, \\ 
&  \\ 
\;\;\;\;\;\;\;\;\;\;\;\;\;\;\;\;\;\;\;\cdots \cdots \cdots \cdots \cdots
\cdots \cdots \cdots \cdots &  \\ 
&  \\ 
t-\left[ \left( p_{1}+\cdots +p_{n-1}\right) b+p_{n}a\right] , & a\leq
x_{n-1}\leq t\leq x_{n}\leq b, \\ 
t-b, & a\leq x_{n}<t\leq b,%
\end{array}%
\right.
\end{equation*}%
\begin{equation*}
=\left\{ 
\begin{array}{ll}
t-a, & a\leq t\leq x_{1}\leq b, \\ 
t-\left( P_{1}b+\bar{P}_{1}a\right) , & a\leq x_{1}<t\leq x_{2}\leq b, \\ 
\cdots \cdots \cdots \cdots \cdots \cdots \cdots \cdots \cdots &  \\ 
t-\left( P_{i}b+\bar{P}_{i}a\right) & a\leq x_{i}\leq t\leq x_{i+1}\leq b,
\\ 
\cdots \cdots \cdots \cdots \cdots \cdots \cdots \cdots \cdots &  \\ 
t-\left( P_{n-1}b+\bar{P}_{n-1}a\right) , & a\leq x_{n-1}\leq t\leq
x_{n}\leq b, \\ 
t-b, & a\leq x_{n}<t\leq b.%
\end{array}%
\right.
\end{equation*}%
Consequently, by (\ref{2.3}) and (\ref{2.4}), we have 
\begin{align}
& \sum_{i=1}^{n}p_{i}f\left( x_{i}\right)  \label{2.5} \\
& =\frac{1}{b-a}\left( B\right) \int_{a}^{b}f\left( t\right) dt+\frac{1}{b-a}%
\left( B\right) \int_{a}^{b}S\left( \bar{x},\bar{p},t\right) f^{\prime
}\left( t\right) dt\;\;\;\text{(by (\ref{2.4}))}  \notag \\
& =\frac{1}{b-a}\left( B\right) \int_{a}^{b}f\left( t\right) dt+\frac{1}{b-a}%
\left( B\right) \int_{a}^{x_{1}}\left( t-a\right) f^{\prime }\left( t\right)
dt  \notag \\
& \qquad \qquad +\frac{1}{b-a}\sum_{i=1}^{n-1}\left( B\right) \int_{a}^{b}%
\left[ t-\left( P_{i}b+\bar{P}_{i}a\right) \right] f^{\prime }\left(
t\right) dt  \notag \\
& \qquad \qquad +\frac{1}{b-a}\left( B\right) \int_{x_{n}}^{b}\left(
t-b\right) f^{\prime }\left( t\right) dt,  \notag
\end{align}%
and the representation (\ref{2.1}) is proved.
\end{proof}

The following result in approximating the Bochner integral $\left( B\right)
\int_{a}^{b}f\left( t\right) dt$ in terms of the convex combination of $%
\left( f\left( x_{i}\right) \right) _{i=\overline{1,n}}$ with the weights $%
\left( p_{i}\right) _{i=\overline{1,n}}$ holds.

\begin{theorem}
\label{t1}Assume that $f:\left[ a,b\right] \rightarrow X,$ $\left(
x_{i}\right) _{i=\overline{1,n}}$ and $\left( p_{i}\right) _{i=\overline{1,n}%
}$ are as in Lemma \ref{l1}. Then we have the inequality: 
\begin{multline}
\left\Vert \left( B\right) \int_{a}^{b}f\left( t\right) dt-\left( b-a\right)
\sum_{i=1}^{n}p_{i}f\left( x_{i}\right) \right\Vert  \label{2.6} \\
\leq \int_{a}^{x_{1}}\left( t-a\right) \left\Vert f^{\prime }\left( t\right)
\right\Vert dt+\sum_{i=1}^{n-1}\int_{x_{i}}^{x_{i+1}}\left\vert t-\left(
P_{i}b+\bar{P}_{i}a\right) \right\vert \left\Vert f^{\prime }\left( t\right)
\right\Vert dt \\
+\int_{x_{n}}^{b}\left( b-t\right) \left\Vert f^{\prime }\left( t\right)
\right\Vert dt
\end{multline}%
\begin{equation*}
\leq \left\{ 
\begin{array}{l}
\left( x_{1}-a\right) \left\Vert \left\vert f^{\prime }\right\vert
\right\Vert _{\left[ a,x_{1}\right] ,1}+\sum\limits_{i=1}^{n-1}\mu _{\infty
}\left( x_{i},P_{i}b+\bar{P}_{i}a,x_{i+1}\right) \left\Vert \left\vert
f^{\prime }\right\vert \right\Vert _{\left[ x_{i},x_{i+1}\right] ,1} \\ 
\hfill +\left( b-x_{n}\right) \left\Vert \left\vert f^{\prime }\right\vert
\right\Vert _{\left[ x_{n},b\right] ,1} \\ 
\\ 
\frac{\left( x_{1}-a\right) ^{1+\frac{1}{q}}}{\left( q+1\right) ^{\frac{1}{q}%
}}\left\Vert \left\vert f^{\prime }\right\vert \right\Vert _{\left[ a,x_{1}%
\right] ,p}+\sum\limits_{i=1}^{n-1}\left[ \mu _{q}\left( x_{i},P_{i}b+\bar{P}%
_{i}a,x_{i+1}\right) \right] ^{\frac{1}{q}}\left\Vert \left\vert f^{\prime
}\right\vert \right\Vert _{\left[ x_{i},x_{i+1}\right] ,p} \\ 
\hfill +\frac{\left( b-x_{n}\right) ^{1+\frac{1}{q}}}{\left( q+1\right) ^{%
\frac{1}{q}}}\left\Vert \left\vert f^{\prime }\right\vert \right\Vert _{%
\left[ x_{n},b\right] ,p}\text{ \hspace{0.05in}where }p>1,\;\frac{1}{p}+%
\frac{1}{q}=1\text{ } \\ 
\hfill \text{and if }f^{\prime }\in L_{p}\left( \left[ a,b\right] ;X\right) ;
\\ 
\\ 
\frac{\left( x_{1}-a\right) ^{2}}{2}\left\Vert \left\vert f^{\prime
}\right\vert \right\Vert _{\left[ a,x_{1}\right] ,\infty
}+\sum\limits_{i=1}^{n-1}\mu _{1}\left( x_{i},P_{i}b+\bar{P}%
_{i}a,x_{i+1}\right) \left\Vert \left\vert f^{\prime }\right\vert
\right\Vert _{\left[ x_{i},x_{i+1}\right] ,\infty } \\ 
\hfill +\frac{\left( b-x_{n}\right) ^{2}}{2}\left\Vert \left\vert f^{\prime
}\right\vert \right\Vert _{\left[ x_{n},b\right] ,\infty }\;\;\text{ if }%
f^{\prime }\in L_{\infty }\left( \left[ a,b\right] ;X\right) ;%
\end{array}%
\right.
\end{equation*}%
\begin{equation*}
\leq \left\{ 
\begin{array}{l}
\max \left( x_{1}-a,\max\limits_{i=\overline{1,n-1}}\left\{ \mu _{\infty
}\left( x_{i},P_{i}b+\bar{P}_{i}a,x_{i+1}\right) \right\} ,b-x_{n}\right)
\left\Vert \left\vert f^{\prime }\right\vert \right\Vert _{\left[ a,b\right]
,1} \\ 
\hfill \text{if }f^{\prime }\in L_{1}\left( \left[ a,b\right] ;X\right) ; \\ 
\\ 
\left[ \frac{\left( x_{1}-a\right) ^{q+1}}{q+1}+\sum\limits_{i=1}^{n-1}\mu
_{q}\left( x_{i},P_{i}b+\bar{P}_{i}a,x_{i+1}\right) +\frac{\left(
b-x_{n}\right) ^{q+1}}{q+1}\right] ^{\frac{1}{q}}\left\Vert \left\vert
f^{\prime }\right\vert \right\Vert _{\left[ a,b\right] ,p} \\ 
\hfill \text{where }p>1,\;\frac{1}{p}+\frac{1}{q}=1\text{ }\;\text{and if }%
f^{\prime }\in L_{p}\left( \left[ a,b\right] ;X\right) ; \\ 
\\ 
\left[ \frac{\left( x_{1}-a\right) ^{2}}{2}+\sum\limits_{i=1}^{n-1}\mu
_{1}\left( x_{i},P_{i}b+\bar{P}_{i}a,x_{i+1}\right) +\frac{\left(
b-x_{n}\right) ^{2}}{2}\right] \left\Vert \left\vert f^{\prime }\right\vert
\right\Vert _{\left[ a,b\right] ,\infty } \\ 
\hfill \text{if }f^{\prime }\in L_{\infty }\left( \left[ a,b\right]
;X\right) ;%
\end{array}%
\right.
\end{equation*}%
where $L_{p}\left( \left[ a,b\right] ;X\right) ,$ $p\in \left[ 1,\infty %
\right] $ are the usual vector-valued Lebesgue spaces and 
\begin{align*}
\left\Vert \left\vert h\right\vert \right\Vert _{\left[ \alpha ,\beta \right]
,\infty }& :=ess\sup\limits_{t\in \left[ \alpha ,\beta \right] }\left\Vert
h\left( t\right) \right\Vert , \\
\left\Vert \left\vert h\right\vert \right\Vert _{\left[ \alpha ,\beta \right]
,p}& :=\left( \int_{\alpha }^{\beta }\left\Vert h\left( t\right) \right\Vert
^{p}dt\right) ^{\frac{1}{p}},\;p\geq 1,
\end{align*}%
and the functions $\mu _{q}\left( \cdot ,\cdot ,\cdot \right) ,$ $q\in \left[
1,\infty \right] $ were defined in (\ref{2.0}).
\end{theorem}

\begin{proof}
Using the properties of the norm, we have, by (\ref{2.1}), that 
\begin{multline*}
\left\| \left( B\right) \int_{a}^{b}f\left( t\right) dt-\left( b-a\right)
\sum_{i=1}^{n}p_{i}f\left( x_{i}\right) \right\| \\
\leq \int_{a}^{x_{1}}\left( t-a\right) \left\| f^{\prime }\left( t\right)
\right\| dt+\sum_{i=1}^{n-1}\int_{x_{i}}^{x_{i+1}}\left| t-\left( P_{i}b+%
\bar{P}_{i}a\right) \right| \left\| f^{\prime }\left( t\right) \right\| dt \\
+\int_{x_{n}}^{b}\left( b-t\right) \left\| f^{\prime }\left( t\right)
\right\| dt
\end{multline*}
and the first inequality in (\ref{2.6}) is proved.

Now, observe that 
\begin{multline*}
\int_{a}^{x_{1}}\left( t-a\right) \left\Vert f^{\prime }\left( t\right)
\right\Vert dt \\
\leq \left\{ 
\begin{array}{l}
\left( x_{1}-a\right) \left\Vert \left\vert f^{\prime }\right\vert
\right\Vert _{\left[ a,x_{1}\right] ,1}; \\ 
\\ 
\frac{\left( x_{1}-a\right) ^{1+\frac{1}{q}}}{\left( q+1\right) ^{\frac{1}{q}%
}}\left\Vert \left\vert f^{\prime }\right\vert \right\Vert _{\left[ a,x_{1}%
\right] ,p}\text{ \hspace{0.05in}if }p>1,\;\frac{1}{p}+\frac{1}{q}=1\text{, }%
f^{\prime }\in L_{p}\left( \left[ a,b\right] ;X\right) ; \\ 
\\ 
\frac{\left( x_{1}-a\right) ^{2}}{2}\left\Vert \left\vert f^{\prime
}\right\vert \right\Vert _{\left[ a,x_{1}\right] ,\infty }\;\;\text{if }%
f^{\prime }\in L_{\infty }\left( \left[ a,b\right] ;X\right) ;%
\end{array}%
\right.
\end{multline*}%
and 
\begin{align*}
& \sum_{i=1}^{n-1}\int_{x_{i}}^{x_{i+1}}\left\vert t-\left( P_{i}b+\bar{P}%
_{i}a\right) \right\vert \left\Vert f^{\prime }\left( t\right) \right\Vert dt
\\
& \\
& \leq \left\{ 
\begin{array}{l}
\sup\limits_{t\in \left[ x_{i},x_{i+1}\right] }\left\vert t-\left( P_{i}b+%
\bar{P}_{i}a\right) \right\vert \left\Vert \left\vert f^{\prime }\right\vert
\right\Vert _{\left[ x_{i},x_{i+1}\right] ,1}; \\ 
\\ 
\left( \int_{x_{i}}^{x_{i+1}}\left\vert t-\left( P_{i}b+\bar{P}_{i}a\right)
\right\vert ^{q}dt\right) ^{\frac{1}{q}}\left\Vert \left\vert f^{\prime
}\right\vert \right\Vert _{\left[ x_{i},x_{i+1}\right] ,p} \\ 
\text{ \hspace{0.05in}if }p>1,\;\frac{1}{p}+\frac{1}{q}=1\text{, and }%
f^{\prime }\in L_{p}\left( \left[ a,b\right] ;X\right) ; \\ 
\\ 
\int_{x_{i}}^{x_{i+1}}\left\vert t-\left( P_{i}b+\bar{P}_{i}a\right)
\right\vert dt\left\Vert \left\vert f^{\prime }\right\vert \right\Vert _{%
\left[ x_{i},x_{i+1}\right] ,\infty }\;\;\text{if }f^{\prime }\in L_{\infty
}\left( \left[ a,b\right] ;X\right) ;%
\end{array}%
\right. \\
& \\
& =\left\{ 
\begin{array}{l}
\mu _{\infty }\left( x_{i},P_{i}b+\bar{P}_{i}a,x_{i+1}\right) \left\Vert
\left\vert f^{\prime }\right\vert \right\Vert _{\left[ x_{i},x_{i+1}\right]
,1}; \\ 
\\ 
\left[ \mu _{q}\left( x_{i},P_{i}b+\bar{P}_{i}a,x_{i+1}\right) \right] ^{%
\frac{1}{q}}\left\Vert \left\vert f^{\prime }\right\vert \right\Vert _{\left[
x_{i},x_{i+1}\right] ,p} \\ 
\text{ \hspace{0.05in}if }p>1,\;\frac{1}{p}+\frac{1}{q}=1\text{, and }%
f^{\prime }\in L_{p}\left( \left[ a,b\right] ;X\right) ; \\ 
\\ 
\mu _{1}\left( x_{i},P_{i}b+\bar{P}_{i}a,x_{i+1}\right) \left\Vert
\left\vert f^{\prime }\right\vert \right\Vert _{\left[ x_{i},x_{i+1}\right]
,\infty }\;\;\text{if }f^{\prime }\in L_{\infty }\left( \left[ a,b\right]
;X\right) ;%
\end{array}%
\right.
\end{align*}%
with 
\begin{multline*}
\int_{x_{n}}^{b}\left( b-t\right) \left\Vert f^{\prime }\left( t\right)
\right\Vert dt \\
\leq \left\{ 
\begin{array}{l}
\left( b-x_{n}\right) \left\Vert \left\vert f^{\prime }\right\vert
\right\Vert _{\left[ x_{n},b\right] ,1}; \\ 
\\ 
\frac{\left( b-x_{n}\right) ^{1+\frac{1}{q}}}{\left( q+1\right) ^{\frac{1}{q}%
}}\left\Vert \left\vert f^{\prime }\right\vert \right\Vert _{\left[ x_{n},b%
\right] ,p}\text{ \hspace{0.05in}if }p>1,\;\frac{1}{p}+\frac{1}{q}=1\text{, }%
f^{\prime }\in L_{p}\left( \left[ a,b\right] ;X\right) ; \\ 
\\ 
\frac{\left( b-x_{n}\right) ^{2}}{2}\left\Vert \left\vert f^{\prime
}\right\vert \right\Vert _{\left[ x_{n},b\right] ,\infty }\;\;\text{if }%
f^{\prime }\in L_{\infty }\left( \left[ a,b\right] ;X\right) ;%
\end{array}%
\right.
\end{multline*}%
giving the second inequality in (\ref{2.6}).

Finally, observe that 
\begin{eqnarray*}
&&\left( x_{1}-a\right) \left\Vert \left\vert f^{\prime }\right\vert
\right\Vert _{\left[ a,x_{1}\right] ,1}+\sum\limits_{i=1}^{n-1}\mu _{\infty
}\left( x_{i},P_{i}b+\bar{P}_{i}a,x_{i+1}\right) \left\Vert \left\vert
f^{\prime }\right\vert \right\Vert _{\left[ x_{i},x_{i+1}\right] ,1} \\
&&+\left( b-x_{n}\right) \left\Vert \left\vert f^{\prime }\right\vert
\right\Vert _{\left[ x_{n},b\right] ,1} \\
&\leq &\max \left\{ x_{1}-a,\max\limits_{i=\overline{1,n-1}}\left\{ \mu
_{\infty }\left( x_{i},P_{i}b+\bar{P}_{i}a,x_{i+1}\right) \right\}
,b-x_{n}\right\} \\
&&\times \left[ \left\Vert \left\vert f^{\prime }\right\vert \right\Vert _{%
\left[ a,x_{1}\right] ,1}+\sum\limits_{i=1}^{n-1}\left\Vert \left\vert
f^{\prime }\right\vert \right\Vert _{\left[ x_{i},x_{i+1}\right]
,1}+\left\Vert \left\vert f^{\prime }\right\vert \right\Vert _{\left[ x_{n},b%
\right] ,1}\right] \\
&\leq &\max \left\{ x_{1}-a,\max\limits_{i=\overline{1,n-1}}\left\{ \mu
_{\infty }\left( x_{i},P_{i}b+\bar{P}_{i}a,x_{i+1}\right) \right\}
,b-x_{n}\right\} \left\Vert \left\vert f^{\prime }\right\vert \right\Vert _{%
\left[ a,b\right] ,1}.
\end{eqnarray*}%
Further, by the discrete H\"{o}lder's inequality we have that 
\begin{eqnarray*}
&&\frac{\left( x_{1}-a\right) ^{1+\frac{1}{q}}}{\left( q+1\right) ^{\frac{1}{%
q}}}\left\Vert \left\vert f^{\prime }\right\vert \right\Vert _{\left[ a,x_{1}%
\right] ,p}+\sum\limits_{i=1}^{n-1}\left[ \mu _{q}\left( x_{i},P_{i}b+\bar{P}%
_{i}a,x_{i+1}\right) \right] ^{\frac{1}{q}}\left\Vert \left\vert f^{\prime
}\right\vert \right\Vert _{\left[ x_{i},x_{i+1}\right] ,p} \\
&&+\frac{\left( b-x_{n}\right) ^{1+\frac{1}{q}}}{\left( q+1\right) ^{\frac{1%
}{q}}}\left\Vert \left\vert f^{\prime }\right\vert \right\Vert _{\left[
x_{n},b\right] ,p} \\
&\leq &\left\{ \left[ \frac{\left( x_{1}-a\right) ^{1+\frac{1}{q}}}{\left(
q+1\right) ^{\frac{1}{q}}}\right] ^{q}+\sum\limits_{i=1}^{n-1}\left( \left[
\mu _{q}\left( x_{i},P_{i}b+\bar{P}_{i}a,x_{i+1}\right) \right] ^{\frac{1}{q}%
}\right) ^{q}\right. \\
&&+\left. \left[ \frac{\left( b-x_{n}\right) ^{1+\frac{1}{q}}}{\left(
q+1\right) ^{\frac{1}{q}}}\right] ^{q}\right\} ^{\frac{1}{q}}\times \left[
\left\Vert \left\vert f^{\prime }\right\vert \right\Vert _{\left[ x_{n},b%
\right] ,p}^{p}+\left\Vert \left\vert f^{\prime }\right\vert \right\Vert _{%
\left[ a,x_{1}\right] ,p}^{p}+\sum\limits_{i=1}^{n-1}\left\Vert \left\vert
f^{\prime }\right\vert \right\Vert _{\left[ x_{i},x_{i+1}\right] ,p}^{p}%
\right] ^{\frac{1}{p}} \\
&=&\left[ \frac{\left( x_{1}-a\right) ^{q+1}}{q+1}+\sum\limits_{i=1}^{n-1}%
\mu _{q}\left( x_{i},P_{i}b+\bar{P}_{i}a,x_{i+1}\right) +\frac{\left(
b-x_{n}\right) ^{q+1}}{q+1}\right] ^{\frac{1}{q}}\left\Vert \left\vert
f^{\prime }\right\vert \right\Vert _{\left[ a,b\right] ,p}
\end{eqnarray*}%
and 
\begin{eqnarray*}
&&\frac{\left( x_{1}-a\right) ^{2}}{2}\left\Vert \left\vert f^{\prime
}\right\vert \right\Vert _{\left[ a,x_{1}\right] ,\infty
}+\sum\limits_{i=1}^{n-1}\mu _{1}\left( x_{i},P_{i}b+\bar{P}%
_{i}a,x_{i+1}\right) \left\Vert \left\vert f^{\prime }\right\vert
\right\Vert _{\left[ x_{i},x_{i+1}\right] ,\infty } \\
&&+\frac{\left( b-x_{n}\right) ^{2}}{2}\left\Vert \left\vert f^{\prime
}\right\vert \right\Vert _{\left[ x_{n},b\right] ,\infty } \\
&\leq &\left[ \frac{\left( x_{1}-a\right) ^{2}}{2}+\sum\limits_{i=1}^{n-1}%
\mu _{1}\left( x_{i},P_{i}b+\bar{P}_{i}a,x_{i+1}\right) +\frac{\left(
b-x_{n}\right) ^{2}}{2}\right] \\
&&\times \max \left\{ \left\Vert \left\vert f^{\prime }\right\vert
\right\Vert _{\left[ a,x_{1}\right] ,\infty },\max\limits_{i=\overline{1,n-1}%
}\left\Vert \left\vert f^{\prime }\right\vert \right\Vert _{\left[
x_{i},x_{i+1}\right] ,\infty },\left\Vert \left\vert f^{\prime }\right\vert
\right\Vert _{\left[ x_{n},b\right] ,\infty }\right\} \\
&=&\left[ \frac{\left( x_{1}-a\right) ^{2}}{2}+\sum\limits_{i=1}^{n-1}\mu
_{1}\left( x_{i},P_{i}b+\bar{P}_{i}a,x_{i+1}\right) +\frac{\left(
b-x_{n}\right) ^{2}}{2}\right] \left\Vert \left\vert f^{\prime }\right\vert
\right\Vert _{\left[ a,b\right] ,\infty };
\end{eqnarray*}%
and the theorem is completely proved.
\end{proof}

It is a natural assumption to consider the weights $p_{i}>0$ $\left(
i=1,\dots ,n\right) $ for which $\xi _{i}:=P_{i}b+\bar{P}_{i}a$ $\left( \in %
\left[ a,b\right] \right) $ will be in the interval $\left[ x_{i},x_{i+1}%
\right] $ \hspace{0.05in}$\left( i=1,\dots ,n\right) .$ In this case we
have: 
\begin{equation*}
\mu _{\infty }\left( x_{i},P_{i}b+\bar{P}_{i}a,x_{i+1}\right) =\frac{1}{2}%
h_{i}+\left| P_{i}b+\bar{P}_{i}a-\frac{x_{i}+x_{i+1}}{2}\right| ,
\end{equation*}
where $h_{i}:=x_{i+1}-x_{i}$, and for $p\in \lbrack 1,\infty )$%
\begin{equation*}
\mu _{p}\left( x_{i},P_{i}b+\bar{P}_{i}a,x_{i+1}\right) =\frac{1}{p+1}\left[
\left( P_{i}b+\bar{P}_{i}a-x_{i}\right) ^{p+1}+\left( x_{i+1}-P_{i}b-\bar{P}%
_{i}a\right) ^{p+1}\right] .
\end{equation*}
Note that for $p=1,$ we have 
\begin{equation*}
\mu _{1}\left( x_{i},P_{i}b+\bar{P}_{i}a,x_{i+1}\right) =\frac{1}{4}%
h_{i}^{2}+\left( P_{i}b+\bar{P}_{i}a-\frac{x_{i}+x_{i+1}}{2}\right) ^{2}.
\end{equation*}

The following corollary is important for applications.

\begin{corollary}
\label{c1}With the assumptions of Lemma \ref{l1} and if $x_{i}\leq P_{i}b+%
\bar{P}_{i}a\leq x_{i+1}$ for each $i=1,\dots ,n-1,$ then we have the
inequalities: 
\begin{equation}
\left\Vert \left( B\right) \int_{a}^{b}f\left( t\right) dt-\left( b-a\right)
\sum_{i=1}^{n}p_{i}f\left( x_{i}\right) \right\Vert  \label{2.7}
\end{equation}%
\begin{eqnarray*}
&\leq &\left\{ 
\begin{array}{l}
\left( x_{1}-a\right) \left\Vert \left\vert f^{\prime }\right\vert
\right\Vert _{\left[ a,x_{1}\right] ,1}+\sum\limits_{i=1}^{n-1}\left[ \frac{1%
}{2}h_{i}+\left\vert P_{i}b+\bar{P}_{i}a-\frac{x_{i}+x_{i+1}}{2}\right\vert %
\right] \left\Vert \left\vert f^{\prime }\right\vert \right\Vert _{\left[
x_{i},x_{i+1}\right] ,1} \\ 
\hfill +\left( b-x_{n}\right) \left\Vert \left\vert f^{\prime }\right\vert
\right\Vert _{\left[ x_{n},b\right] ,1}; \\ 
\\ 
\frac{\left( x_{1}-a\right) ^{1+\frac{1}{q}}}{\left( q+1\right) ^{\frac{1}{q}%
}}\left\Vert \left\vert f^{\prime }\right\vert \right\Vert _{\left[ a,x_{1}%
\right] ,p} \\ 
+\frac{1}{\left( q+1\right) ^{\frac{1}{q}}}\sum\limits_{i=1}^{n-1}\left[
\left( P_{i}b+\bar{P}_{i}a-x_{i}\right) ^{q+1}+\left( x_{i+1}-P_{i}b-\bar{P}%
_{i}a\right) ^{q+1}\right] ^{\frac{1}{q}}\left\Vert \left\vert f^{\prime
}\right\vert \right\Vert _{\left[ x_{i},x_{i+1}\right] ,p} \\ 
\hfill +\frac{\left( b-x_{n}\right) ^{1+\frac{1}{q}}}{\left( q+1\right) ^{%
\frac{1}{q}}}\left\Vert \left\vert f^{\prime }\right\vert \right\Vert _{%
\left[ x_{n},b\right] ,p}\text{ \hspace{0.05in}where }p>1,\;\frac{1}{p}+%
\frac{1}{q}=1\text{ } \\ 
\hfill \text{and if }f^{\prime }\in L_{p}\left( \left[ a,b\right] ;X\right) ;
\\ 
\\ 
\frac{\left( x_{1}-a\right) ^{2}}{2}\left\Vert \left\vert f^{\prime
}\right\vert \right\Vert _{\left[ a,x_{1}\right] ,\infty
}+\sum\limits_{i=1}^{n-1}\left[ \frac{1}{4}h_{i}^{2}+\left( P_{i}b+\bar{P}%
_{i}a-\frac{x_{i}+x_{i+1}}{2}\right) ^{2}\right] \left\Vert \left\vert
f^{\prime }\right\vert \right\Vert _{\left[ x_{i},x_{i+1}\right] ,\infty }
\\ 
\hfill +\frac{\left( b-x_{n}\right) ^{2}}{2}\left\Vert \left\vert f^{\prime
}\right\vert \right\Vert _{\left[ x_{n},b\right] ,\infty }\;\;\text{ if }%
f^{\prime }\in L_{\infty }\left( \left[ a,b\right] ;X\right) ;%
\end{array}%
\right. \\
&&
\end{eqnarray*}%
\begin{equation*}
\leq \left\{ 
\begin{array}{l}
\max \left( x_{1}-a,\frac{1}{2}\max\limits_{i=\overline{1,n-1}%
}h_{i}+\max\limits_{i=\overline{1,n-1}}\left\vert P_{i}b+\bar{P}_{i}a-\frac{%
x_{i}+x_{i+1}}{2}\right\vert ,b-x_{n}\right) \left\Vert \left\vert f^{\prime
}\right\vert \right\Vert _{\left[ a,b\right] ,1} \\ 
\hfill \text{if }f^{\prime }\in L_{1}\left( \left[ a,b\right] ;X\right) ; \\ 
\\ 
\frac{1}{\left( q+1\right) ^{\frac{1}{q}}}\left[ \left( x_{1}-a\right)
^{q+1}+\sum\limits_{i=1}^{n-1}\right. \left[ \left( P_{i}b+\bar{P}%
_{i}a-x_{i}\right) ^{q+1}\right. \\ 
\hfill +\left. \left( x_{i+1}-P_{i}b-\bar{P}_{i}a\right) ^{q+1}\right]
+\left( b-x_{n}\right) ^{q+1}\bigg]^{\frac{1}{q}}\left\Vert \left\vert
f^{\prime }\right\vert \right\Vert _{\left[ a,b\right] ,p} \\ 
\hfill \text{if }p>1,\;\frac{1}{p}+\frac{1}{q}=1\text{,\ }f^{\prime }\in
L_{p}\left( \left[ a,b\right] ;X\right) ; \\ 
\\ 
\left[ \frac{\left( x_{1}-a\right) ^{2}}{2}+\sum\limits_{i=1}^{n-1}\left[ 
\frac{1}{4}h_{i}^{2}+\left( P_{i}b+\bar{P}_{i}a-\frac{x_{i}+x_{i+1}}{2}%
\right) ^{2}\right] +\frac{\left( b-x_{n}\right) ^{2}}{2}\right] \left\Vert
\left\vert f^{\prime }\right\vert \right\Vert _{\left[ a,b\right] ,\infty }
\\ 
\hfill \text{if }f^{\prime }\in L_{\infty }\left( \left[ a,b\right]
;X\right) .%
\end{array}%
\right.
\end{equation*}
\end{corollary}

\begin{remark}
For $n=1,$ we recapture from $\left( \ref{2.7}\right) $ the Ostrowski type
inequalities incorporated in Theorems \ref{t2.1} and \ref{ta}.
\end{remark}

\section{The Case of Two Points}

The following proposition is a particular case of Corollary \ref{c1} for $%
n=2 $ and will be considered in some details since there are important for
applications.

\begin{proposition}
Let $\left( X,\left\Vert \cdot \right\Vert \right) $ be a Banach space with
the Radon-Nikodym property and $f:\left[ a,b\right] \rightarrow X$ be an
absolutely continuous function on $\left[ a,b\right] .$ If $a\leq x_{1}\leq
x_{2}\leq b$ $\left( b>a\right) $ and $t\in \left[ 0,1\right] $ satisfies
the condition%
\begin{equation*}
\left( 0\leq \right) \frac{x_{1}-a}{b-a}\leq t\leq \frac{x_{2}-a}{b-a}\left(
\leq 1\right) ,
\end{equation*}%
then we have the inequalities%
\begin{equation}
\left\Vert \left( B\right) \int_{a}^{b}f\left( s\right) ds-\left( b-a\right) 
\left[ tf\left( x_{2}\right) +\left( 1-t\right) f\left( x_{1}\right) \right]
\right\Vert  \label{3.2}
\end{equation}%
\begin{equation*}
\leq \left\{ 
\begin{array}{l}
\left( x_{1}-a\right) \left\Vert \left\vert f^{\prime }\right\vert
\right\Vert _{\left[ a,x_{1}\right] ,1}+\left[ \frac{1}{2}\left(
x_{2}-x_{1}\right) +\left\vert tb+\left( 1-t\right) a-\frac{x_{1}+x_{2}}{2}%
\right\vert \right] \left\Vert \left\vert f^{\prime }\right\vert \right\Vert
_{\left[ x_{1},x_{2}\right] ,1} \\ 
+\left( b-x_{2}\right) \left\Vert \left\vert f^{\prime }\right\vert
\right\Vert _{\left[ x_{2},b\right] ,1}; \\ 
\\ 
\frac{1}{\left( q+1\right) ^{1/q}}\left\{ \left( x_{1}-a\right)
^{1+1/q}\left\Vert \left\vert f^{\prime }\right\vert \right\Vert _{\left[
a,x_{1}\right] ,p}\right. \\ 
+\left[ \left( tb+\left( 1-t\right) a-x_{1}\right) ^{q+1}+\left(
x_{2}-tb-\left( 1-t\right) a\right) ^{q+1}\right] ^{1/q}\left\Vert
\left\vert f^{\prime }\right\vert \right\Vert _{\left[ x_{1},x_{2}\right] ,p}
\\ 
\left. +\left( b-x_{2}\right) ^{1+1/q}\left\Vert \left\vert f^{\prime
}\right\vert \right\Vert _{\left[ x_{2},b\right] ,p}\right\} \text{, }p>1,%
\frac{1}{p}+\frac{1}{q}=1,f^{\prime }\in L_{p}\left( \left[ a,b\right]
;X\right) ; \\ 
\\ 
\frac{\left( x_{1}-a\right) ^{2}}{2}\left\Vert \left\vert f^{\prime
}\right\vert \right\Vert _{\left[ a,x_{1}\right] ,\infty }+\left[ \frac{1}{4}%
\left( x_{2}-x_{1}\right) ^{2}+\left[ tb+\left( 1-t\right) a-\frac{%
x_{1}+x_{2}}{2}\right] ^{2}\right] \left\Vert \left\vert f^{\prime
}\right\vert \right\Vert _{\left[ x_{1},x_{2}\right] ,\infty } \\ 
+\frac{\left( b-x_{2}\right) ^{2}}{2}\left\Vert \left\vert f^{\prime
}\right\vert \right\Vert _{\left[ x_{2},b\right] ,\infty }\text{, }f^{\prime
}\in L_{\infty }\left( \left[ a,b\right] ;X\right) ;%
\end{array}%
\right.
\end{equation*}%
$\leq \left\{ 
\begin{array}{l}
\max \left\{ x_{1}-a,\frac{1}{2}\left( x_{2}-x_{1}\right) +\left\vert
tb+\left( 1-t\right) a-\frac{x_{1}+x_{2}}{2}\right\vert ,b-x_{2}\right\}
\left\Vert \left\vert f^{\prime }\right\vert \right\Vert _{\left[ a,b\right]
,1}; \\ 
\\ 
\frac{1}{\left( q+1\right) ^{1/q}}\left\{ \left( x_{1}-a\right)
^{q+1}+\left( tb+\left( 1-t\right) a-x_{1}\right) ^{q+1}+\left(
x_{2}-tb-\left( 1-t\right) a\right) ^{q+1}\right. \\ 
\left. +\left( b-x_{2}\right) ^{q+1}\right\} ^{1/q}\left\Vert \left\vert
f^{\prime }\right\vert \right\Vert _{\left[ a,b\right] ,p}\text{, }p>1,\frac{%
1}{p}+\frac{1}{q}=1,f^{\prime }\in L_{p}\left( \left[ a,b\right] ;X\right) ;
\\ 
\\ 
\left[ \frac{\left( x_{1}-a\right) ^{2}}{2}+\frac{1}{4}\left(
x_{2}-x_{1}\right) ^{2}+\left[ tb+\left( 1-t\right) a-\frac{x_{1}+x_{2}}{2}%
\right] ^{2}+\frac{\left( b-x_{2}\right) ^{2}}{2}\right] \left\Vert
\left\vert f^{\prime }\right\vert \right\Vert _{\left[ a,b\right] ,\infty },
\\ 
\text{ \ \ }f^{\prime }\in L_{\infty }\left( \left[ a,b\right] ;X\right) .%
\end{array}%
\right. $
\end{proposition}

The following particular inequalities are of interest.

\textbf{1}. If $x_{1}=a,x_{2}=b,$ then for any $t\in \left[ 0,1\right] ,$ we
have the inequalities%
\begin{equation}
\left\Vert \left( B\right) \int_{a}^{b}f\left( s\right) ds-\left( b-a\right) 
\left[ tf\left( b\right) +\left( 1-t\right) f\left( a\right) \right]
\right\Vert  \label{3.3}
\end{equation}%
\begin{equation*}
\leq \left\{ 
\begin{array}{l}
\left[ \frac{1}{2}\left( b-a\right) +\left\vert tb+\left( 1-t\right) a-\frac{%
a+b}{2}\right\vert \right] \left\Vert \left\vert f^{\prime }\right\vert
\right\Vert _{\left[ a,b\right] ,1}; \\ 
\\ 
\frac{1}{\left( q+1\right) ^{1/q}}\left[ t^{q+1}+\left( 1-t\right) ^{q+1}%
\right] ^{1/q}\left( b-a\right) ^{1+1/q}\left\Vert \left\vert f^{\prime
}\right\vert \right\Vert _{\left[ a,b\right] ,p}, \\ 
\text{ \ \ }p>1,\frac{1}{p}+\frac{1}{q}=1,f^{\prime }\in L_{p}\left( \left[
a,b\right] ;X\right) ; \\ 
\\ 
\left[ \frac{1}{4}\left( b-a\right) ^{2}+\left( tb+\left( 1-t\right) a-\frac{%
a+b}{2}\right) ^{2}\right] \left\Vert \left\vert f^{\prime }\right\vert
\right\Vert _{\left[ a,b\right] ,\infty },f^{\prime }\in L_{\infty }\left( %
\left[ a,b\right] ;X\right) .%
\end{array}%
\right.
\end{equation*}%
The best inequality one can get from $\left( \ref{3.3}\right) $ is for $t=%
\frac{1}{2},$ obtaining the trapezoidal rule%
\begin{eqnarray}
&&\left\Vert \left( B\right) \int_{a}^{b}f\left( s\right) ds-\left(
b-a\right) .\frac{f\left( b\right) +f\left( a\right) }{2}\right\Vert
\label{3.4} \\
&\leq &\left\{ 
\begin{array}{l}
\frac{1}{2}\left( b-a\right) \left\Vert \left\vert f^{\prime }\right\vert
\right\Vert _{\left[ a,b\right] ,1}; \\ 
\\ 
\frac{1}{2\left( q+1\right) ^{1/q}}\left( b-a\right) ^{1+1/q}\left\Vert
\left\vert f^{\prime }\right\vert \right\Vert _{\left[ a,b\right] ,p},p>1,%
\frac{1}{p}+\frac{1}{q}=1,f^{\prime }\in L_{p}\left( \left[ a,b\right]
;X\right) ; \\ 
\\ 
\frac{1}{4}\left( b-a\right) ^{2}\left\Vert \left\vert f^{\prime
}\right\vert \right\Vert _{\left[ a,b\right] ,\infty },f^{\prime }\in
L_{\infty }\left( \left[ a,b\right] ;X\right) .%
\end{array}%
\right.  \notag
\end{eqnarray}

\textbf{2}. If $x_{1}=\frac{3a+b}{4},x_{2}=\frac{a+3b}{4},$ then for any $%
t\in \left[ \frac{1}{4},\frac{3}{4}\right] $ we have the inequalities%
\begin{equation}
\left\Vert \left( B\right) \int_{a}^{b}f\left( s\right) ds-\left( b-a\right) 
\left[ tf\left( \frac{3a+b}{4}\right) +\left( 1-t\right) f\left( \frac{a+3b}{%
4}\right) \right] \right\Vert  \label{3.5}
\end{equation}%
\begin{equation*}
\leq \left\{ 
\begin{array}{l}
\frac{b-a}{4}\left\Vert \left\vert f^{\prime }\right\vert \right\Vert _{%
\left[ a,\frac{3a+b}{4}\right] ,1}+\left[ \frac{1}{4}\left( b-a\right)
+\left\vert tb+\left( 1-t\right) a-\frac{a+b}{2}\right\vert \right]
\left\Vert \left\vert f^{\prime }\right\vert \right\Vert _{\left[ \frac{3a+b%
}{4},\frac{a+3b}{4}\right] ,1} \\ 
+\frac{b-a}{4}\left\Vert \left\vert f^{\prime }\right\vert \right\Vert _{%
\left[ \frac{a+3b}{4},b\right] ,1}; \\ 
\\ 
\frac{1}{4^{q+1}\left( q+1\right) ^{1/q}}\left( b-a\right)
^{1+1/q}\left\Vert \left\vert f^{\prime }\right\vert \right\Vert _{\left[ a,%
\frac{3a+b}{4}\right] ,p} \\ 
+\frac{1}{\left( q+1\right) ^{1/q}}\left[ \left( tb+\left( 1-t\right) a-%
\frac{3a+b}{4}\right) ^{q+1}\right. \\ 
\left. +\left( \frac{a+3b}{4}-tb-\left( 1-t\right) a\right) ^{q+1}\right]
^{1/q}\left\Vert \left\vert f^{\prime }\right\vert \right\Vert _{\left[ 
\frac{3a+b}{4},\frac{a+3b}{4}\right] ,p} \\ 
+\frac{1}{4^{q+1}\left( q+1\right) ^{1/q}}\left( b-a\right)
^{1+1/q}\left\Vert \left\vert f^{\prime }\right\vert \right\Vert _{\left[ 
\frac{a+3b}{4},b\right] ,p}\text{, }p>1,\frac{1}{p}+\frac{1}{q}=1,f^{\prime
}\in L_{p}\left( \left[ a,b\right] ;X\right) ; \\ 
\\ 
\frac{\left( b-a\right) ^{2}}{8}\left\Vert \left\vert f^{\prime }\right\vert
\right\Vert _{\left[ a,\frac{3a+b}{4}\right] ,\infty }+\left[ \frac{1}{16}%
\left( b-a\right) ^{2}+\left[ tb+\left( 1-t\right) a-\frac{a+b}{2}\right]
^{2}\right] \left\Vert \left\vert f^{\prime }\right\vert \right\Vert _{\left[
\frac{3a+b}{4},\frac{a+3b}{4}\right] ,\infty } \\ 
+\frac{\left( b-a\right) ^{2}}{8}\left\Vert \left\vert f^{\prime
}\right\vert \right\Vert _{\left[ \frac{a+3b}{4},b\right] ,\infty }\text{, }%
f^{\prime }\in L_{\infty }\left( \left[ a,b\right] ;X\right) ;%
\end{array}%
\right.
\end{equation*}%
\begin{equation*}
\leq \left\{ 
\begin{array}{l}
\left[ \frac{1}{4}\left( b-a\right) +\left\vert tb+\left( 1-t\right) a-\frac{%
a+b}{2}\right\vert \right] \left\Vert \left\vert f^{\prime }\right\vert
\right\Vert _{\left[ a,b\right] ,1}; \\ 
\\ 
\frac{1}{\left( q+1\right) ^{1/q}}\left[ \frac{2\left( b-a\right) ^{q+1}}{%
4^{q+1}}+\left( tb+\left( 1-t\right) a-\frac{3a+b}{4}\right) ^{q+1}+\left( 
\frac{a+3b}{4}-tb-\left( 1-t\right) a\right) ^{q+1}\right] ^{1/q} \\ 
\times \left( b-a\right) ^{1+1/q}\left\Vert \left\vert f^{\prime
}\right\vert \right\Vert _{\left[ a,b\right] ,p},p>1,\frac{1}{p}+\frac{1}{q}%
=1,f^{\prime }\in L_{p}\left( \left[ a,b\right] ;X\right) ; \\ 
\\ 
\left[ \frac{1}{8}\left( b-a\right) ^{2}+\left( tb+\left( 1-t\right) a-\frac{%
a+b}{2}\right) ^{2}\right] \left\Vert \left\vert f^{\prime }\right\vert
\right\Vert _{\left[ a,b\right] ,\infty },f^{\prime }\in L_{\infty }\left( %
\left[ a,b\right] ;X\right) .%
\end{array}%
\right.
\end{equation*}%
The best inequality one can get from $\left( \ref{3.5}\right) $ is for $t=%
\frac{1}{2},$ obtaining%
\begin{eqnarray}
&&\left\Vert \left( B\right) \int_{a}^{b}f\left( s\right) ds-\left(
b-a\right) .\frac{f\left( \frac{a+3b}{4}\right) +f\left( \frac{3a+b}{4}%
\right) }{2}\right\Vert  \label{3.6} \\
&\leq &\left\{ 
\begin{array}{l}
\frac{1}{4}\left( b-a\right) \left\Vert \left\vert f^{\prime }\right\vert
\right\Vert _{\left[ a,b\right] ,1}; \\ 
\\ 
\frac{1}{4\left( q+1\right) ^{1/q}}\left( b-a\right) ^{1+1/q}\left\Vert
\left\vert f^{\prime }\right\vert \right\Vert _{\left[ a,b\right] ,p},p>1,%
\frac{1}{p}+\frac{1}{q}=1,f^{\prime }\in L_{p}\left( \left[ a,b\right]
;X\right) ; \\ 
\\ 
\frac{1}{8}\left( b-a\right) ^{2}\left\Vert \left\vert f^{\prime
}\right\vert \right\Vert _{\left[ a,b\right] ,\infty },f^{\prime }\in
L_{\infty }\left( \left[ a,b\right] ;X\right) .%
\end{array}%
\right.  \notag
\end{eqnarray}

\begin{remark}
One may realize that, instead of using the trapezoidal rule in approximating
the Bochner integral $\left( B\right) \int_{a}^{b}f\left( t\right) dt,$ that
one should use the rule%
\begin{equation}
QT\left( f;a,b\right) :=\left( b-a\right) .\frac{f\left( \frac{a+3b}{4}%
\right) +f\left( \frac{3a+b}{4}\right) }{2},  \label{3.7}
\end{equation}%
which provides a halving of the bound on the error.
\end{remark}

\section{The Case of Three Points}

The case of three points is important for applications since it contains
amongst others Simpson's quadrature rule.

The following proposition holds:

\begin{proposition}
\label{4.1}Let $\left( X,\left\Vert \cdot \right\Vert \right) $ be a Banach
space with the Radon-Nikodym property and $f:\left[ a,b\right] \rightarrow X$
be an absolutely continuous function on $\left[ a,b\right] .$ If $a\leq
x_{1}\leq x_{2}\leq x_{3}\leq b$ $\left( b>a\right) $ and $\alpha ,\beta
,\gamma \in \left[ 0,1\right] $ with $\alpha +\beta +\gamma =1$ satisfies
the condition%
\begin{equation}
\left( 0\leq \right) \frac{x_{1}-a}{b-a}\leq \alpha \leq \frac{x_{2}-a}{b-a}%
\leq \alpha +\beta \leq \frac{x_{3}-a}{b-a}\left( \leq 1\right) ,
\label{e4.0}
\end{equation}%
then we have the inequalities%
\begin{equation}
\left\Vert \left( B\right) \int_{a}^{b}f\left( t\right) dt-\left( b-a\right) 
\left[ \alpha f\left( x_{1}\right) +\beta f\left( x_{2}\right) +\left(
1-\alpha -\beta \right) f\left( x_{3}\right) \right] \right\Vert
\label{e4.1}
\end{equation}%
\begin{equation*}
\leq \left\{ 
\begin{array}{l}
\left( x_{1}-a\right) \left\Vert \left\vert f^{\prime }\right\vert
\right\Vert _{\left[ a,x_{1}\right] ,1}+\left[ \frac{1}{2}\left(
x_{2}-x_{1}\right) +\left\vert \alpha b+\left( 1-\alpha \right) a-\frac{%
x_{1}+x_{2}}{2}\right\vert \right] \left\Vert \left\vert f^{\prime
}\right\vert \right\Vert _{\left[ x_{1},x_{2}\right] ,1} \\ 
+\left[ \frac{1}{2}\left( x_{3}-x_{2}\right) +\left\vert \left( \alpha
+\beta \right) b+\left( 1-\alpha -\beta \right) a-\frac{x_{2}+x_{3}}{2}%
\right\vert \right] \left\Vert \left\vert f^{\prime }\right\vert \right\Vert
_{\left[ x_{2},x_{3}\right] ,1} \\ 
+\left( b-x_{2}\right) \left\Vert \left\vert f^{\prime }\right\vert
\right\Vert _{\left[ x_{2},b\right] ,1}; \\ 
\\ 
\frac{1}{\left( q+1\right) ^{1/q}}\left\{ \left( x_{1}-a\right)
^{1+1/q}\left\Vert \left\vert f^{\prime }\right\vert \right\Vert _{\left[
a,x_{1}\right] ,p}\right. \\ 
+\left[ \left( \alpha b+\left( 1-\alpha \right) a-x_{1}\right) ^{q+1}+\left(
x_{2}-\alpha b-\left( 1-\alpha \right) a\right) ^{q+1}\right]
^{1/q}\left\Vert \left\vert f^{\prime }\right\vert \right\Vert _{\left[
x_{1},x_{2}\right] ,p} \\ 
+\left[ \left( \left( \alpha +\beta \right) b+\left( 1-\alpha -\beta \right)
a-x_{2}\right) ^{q+1}+\left( x_{3}-\left( \alpha +\beta \right) b-\left(
1-\alpha -\beta \right) a\right) ^{q+1}\right] ^{1/q}\left\Vert \left\vert
f^{\prime }\right\vert \right\Vert _{\left[ x_{2},x_{3}\right] ,p} \\ 
\left. +\left( b-x_{3}\right) ^{1+1/q}\left\Vert \left\vert f^{\prime
}\right\vert \right\Vert _{\left[ x_{3},b\right] ,p}\right\} \text{, }p>1,%
\frac{1}{p}+\frac{1}{q}=1,f^{\prime }\in L_{p}\left( \left[ a,b\right]
;X\right) ; \\ 
\\ 
\frac{\left( x_{1}-a\right) ^{2}}{2}\left\Vert \left\vert f^{\prime
}\right\vert \right\Vert _{\left[ a,x_{1}\right] ,\infty }+\left[ \frac{1}{4}%
\left( x_{2}-x_{1}\right) ^{2}+\left[ \alpha b+\left( 1-\alpha \right) a-%
\frac{x_{1}+x_{2}}{2}\right] ^{2}\right] \left\Vert \left\vert f^{\prime
}\right\vert \right\Vert _{\left[ x_{1},x_{2}\right] ,\infty } \\ 
\left[ \frac{1}{4}\left( x_{3}-x_{2}\right) ^{2}+\left[ \left( \alpha +\beta
\right) b+\left( 1-\alpha -\beta \right) a-\frac{x_{2}+x_{3}}{2}\right] ^{2}%
\right] \left\Vert \left\vert f^{\prime }\right\vert \right\Vert _{\left[
x_{2},x_{3}\right] ,\infty } \\ 
+\frac{\left( b-x_{3}\right) ^{2}}{2}\left\Vert \left\vert f^{\prime
}\right\vert \right\Vert _{\left[ x_{3},b\right] ,\infty }\text{, }f^{\prime
}\in L_{\infty }\left( \left[ a,b\right] ;X\right) ;%
\end{array}%
\right.
\end{equation*}%
$\leq \left\{ 
\begin{array}{l}
\max \left\{ x_{1}-a,\frac{1}{2}\left( x_{2}-x_{1}\right) +\left\vert \alpha
b+\left( 1-\alpha \right) a-\frac{x_{1}+x_{2}}{2}\right\vert ,\right. \\ 
\left. \frac{1}{2}\left( x_{2}-x_{1}\right) +\left\vert \left( \alpha +\beta
\right) b+\left( 1-\alpha -\beta \right) a-\frac{x_{2}+x_{3}}{2}\right\vert
,b-x_{2}\right\} \left\Vert \left\vert f^{\prime }\right\vert \right\Vert _{%
\left[ a,b\right] ,1}; \\ 
\\ 
\frac{1}{\left( q+1\right) ^{1/q}}\left\{ \left( x_{1}-a\right)
^{q+1}+\left( \alpha b+\left( 1-\alpha \right) a-x_{1}\right) ^{q+1}+\left(
x_{2}-\alpha b-\left( 1-\alpha \right) a\right) ^{q+1}\right. \\ 
+\left( \left( \alpha +\beta \right) b+\left( 1-\alpha -\beta \right)
a-x_{2}\right) ^{q+1}+\left( x_{3}-\left( \alpha +\beta \right) b-\left(
1-\alpha -\beta \right) a\right) ^{q+1} \\ 
\left. +\left( b-x_{3}\right) ^{q+1}\right\} \left\Vert \left\vert f^{\prime
}\right\vert \right\Vert _{\left[ a,b\right] ,p}\text{, }p>1,\frac{1}{p}+%
\frac{1}{q}=1,f^{\prime }\in L_{p}\left( \left[ a,b\right] ;X\right) ; \\ 
\\ 
\left[ \frac{\left( x_{1}-a\right) ^{2}}{2}+\frac{1}{4}\left(
x_{2}-x_{1}\right) ^{2}+\left[ \alpha b+\left( 1-\alpha \right) a-\frac{%
x_{1}+x_{2}}{2}\right] ^{2}\right. \\ 
\left. +\frac{1}{4}\left( x_{3}-x_{2}\right) ^{2}+\left[ \left( \alpha
+\beta \right) b+\left( 1-\alpha -\beta \right) a-\frac{x_{2}+x_{3}}{2}%
\right] ^{2}+\frac{\left( b-x_{3}\right) ^{2}}{2}\right] \left\Vert
\left\vert f^{\prime }\right\vert \right\Vert _{\left[ a,b\right] ,\infty },
\\ 
f^{\prime }\in L_{\infty }\left( \left[ a,b\right] ;X\right) .%
\end{array}%
\right. $
\end{proposition}

The following particular inequalities are of interest.

\textbf{1}. Assume that $x_{1}=a,x_{2}=\frac{a+b}{2},x_{3}=b$ and $\alpha
,\beta \in \left[ 0,1\right] $ so that $0\leq \alpha \leq \frac{1}{2}\leq
\alpha +\beta \leq 1,$ then we have the inequalities

\begin{equation}
\left\Vert \left( B\right) \int_{a}^{b}f\left( t\right) dt-\left( b-a\right) 
\left[ \alpha f\left( a\right) +\beta f\left( \frac{a+b}{2}\right) +\left(
1-\alpha -\beta \right) f\left( b\right) \right] \right\Vert  \label{e4.3}
\end{equation}%
\begin{equation*}
\leq \left\{ 
\begin{array}{l}
\left[ \frac{1}{4}\left( b-a\right) +\left\vert \alpha b+\left( 1-\alpha
\right) a-\frac{3a+b}{4}\right\vert \right] \left\Vert \left\vert f^{\prime
}\right\vert \right\Vert _{\left[ a,\frac{a+b}{2}\right] ,1} \\ 
+\left[ \frac{1}{4}\left( b-a\right) +\left\vert \left( \alpha +\beta
\right) b+\left( 1-\alpha -\beta \right) a-\frac{a+3b}{4}\right\vert \right]
\left\Vert \left\vert f^{\prime }\right\vert \right\Vert _{\left[ \frac{a+b}{%
2},b\right] ,1}; \\ 
\\ 
\\ 
\frac{1}{\left( q+1\right) ^{1/q}}\left\{ \left[ \alpha ^{q+1}\left(
b-a\right) ^{q+1}+\left( \frac{a+b}{2}-\alpha b-\left( 1-\alpha \right)
a\right) ^{q+1}\right] ^{1/q}\left\Vert \left\vert f^{\prime }\right\vert
\right\Vert _{\left[ a,\frac{a+b}{2}\right] ,p}\right. \\ 
+\left. \left[ \left( \left( \alpha +\beta \right) b+\left( 1-\alpha -\beta
\right) a-\frac{a+b}{2}\right) ^{q+1}+\left( 1-\alpha -\beta \right)
^{q+1}\left( b-a\right) ^{q+1}\right] ^{1/q}\right\} \left\Vert \left\vert
f^{\prime }\right\vert \right\Vert _{\left[ \frac{a+b}{2},b\right] ,p} \\ 
\text{ \ \ }p>1,\frac{1}{p}+\frac{1}{q}=1,f^{\prime }\in L_{p}\left( \left[
a,b\right] ;X\right) ; \\ 
\\ 
\\ 
\left[ \frac{1}{16}\left( b-a\right) ^{2}+\left[ \alpha b+\left( 1-\alpha
\right) a-\frac{3a+b}{4}\right] ^{2}\right] \left\Vert \left\vert f^{\prime
}\right\vert \right\Vert _{\left[ a,\frac{a+b}{2}\right] ,\infty } \\ 
\left[ \frac{1}{16}\left( b-a\right) ^{2}+\left[ \left( \alpha +\beta
\right) b+\left( 1-\alpha -\beta \right) a-\frac{a+3b}{4}\right] ^{2}\right]
\left\Vert \left\vert f^{\prime }\right\vert \right\Vert _{\left[ \frac{a+b}{%
2},b\right] ,\infty } \\ 
\text{ \ \ }f^{\prime }\in L_{\infty }\left( \left[ a,b\right] ;X\right) ;%
\end{array}%
\right.
\end{equation*}

\begin{equation*}
\leq \left\{ 
\begin{array}{l}
\max \left\{ \left[ \frac{1}{4}\left( b-a\right) +\left\vert \alpha b+\left(
1-\alpha \right) a-\frac{3a+b}{4}\right\vert \right] ,\right. \\ 
\left. \left[ \frac{1}{4}\left( b-a\right) +\left\vert \left( \alpha +\beta
\right) b+\left( 1-\alpha -\beta \right) a-\frac{a+3b}{4}\right\vert \right]
\right\} \left\Vert \left\vert f^{\prime }\right\vert \right\Vert _{\left[
a,b\right] ,1} \\ 
\\ 
\\ 
\frac{1}{\left( q+1\right) ^{1/q}}\left\{ \alpha ^{q+1}\left( b-a\right)
^{q+1}+\left( \frac{a+b}{2}-\alpha b-\left( 1-\alpha \right) a\right)
^{q+1}\right. \\ 
+\left. \left( \left( \alpha +\beta \right) b+\left( 1-\alpha -\beta \right)
a-\frac{a+b}{2}\right) ^{q+1}+\left( 1-\alpha -\beta \right) ^{q+1}\left(
b-a\right) ^{q+1}\right\} ^{1/q}\left\Vert \left\vert f^{\prime }\right\vert
\right\Vert _{\left[ a,b\right] ,p} \\ 
\text{ \ \ }p>1,\frac{1}{p}+\frac{1}{q}=1,f^{\prime }\in L_{p}\left( \left[
a,b\right] ;X\right) ; \\ 
\\ 
\\ 
\left\{ \frac{1}{8}\left( b-a\right) ^{2}+\left[ \alpha b+\left( 1-\alpha
\right) a-\frac{3a+b}{4}\right] ^{2}\right. \\ 
\left. +\left[ \left( \alpha +\beta \right) b+\left( 1-\alpha -\beta \right)
a-\frac{a+3b}{4}\right] ^{2}\right\} \left\Vert \left\vert f^{\prime
}\right\vert \right\Vert _{\left[ a,b\right] ,\infty } \\ 
\text{ \ \ }f^{\prime }\in L_{\infty }\left( \left[ a,b\right] ;X\right) .%
\end{array}%
\right.
\end{equation*}%
It is easy to see that, the best inequality one can derive from (\ref{e4.3})
is the one for $\alpha =\frac{1}{4}$ and $\beta =\frac{3}{4},$ getting%
\begin{equation}
\left\Vert \left( B\right) \int_{a}^{b}f\left( t\right) dt-\left( b-a\right) 
\left[ \frac{f\left( a\right) +f\left( b\right) }{4}+\frac{1}{2}f\left( 
\frac{a+b}{2}\right) \right] \right\Vert  \label{e4.4}
\end{equation}%
\begin{equation*}
\leq \left\{ 
\begin{array}{l}
\frac{1}{4}\left( b-a\right) \left\Vert \left\vert f^{\prime }\right\vert
\right\Vert _{\left[ a,b\right] ,1}; \\ 
\\ 
\frac{1}{2^{2+1/q}\left( q+1\right) ^{1/q}}\left( b-a\right) ^{1+1/q}\left\{
\left\Vert \left\vert f^{\prime }\right\vert \right\Vert _{\left[ a,\frac{a+b%
}{2}\right] ,p}+\left\Vert \left\vert f^{\prime }\right\vert \right\Vert _{%
\left[ \frac{a+b}{2},b\right] ,p}\right\} \\ 
\text{ \ \ }p>1,\frac{1}{p}+\frac{1}{q}=1,f^{\prime }\in L_{p}\left( \left[
a,b\right] ;X\right) ; \\ 
\\ 
\\ 
\frac{1}{16}\left( b-a\right) ^{2}\left[ \left\Vert \left\vert f^{\prime
}\right\vert \right\Vert _{\left[ a,\frac{a+b}{2}\right] ,\infty
}+\left\Vert \left\vert f^{\prime }\right\vert \right\Vert _{\left[ \frac{a+b%
}{2},b\right] ,\infty }\right] ,\text{ }f^{\prime }\in L_{\infty }\left( %
\left[ a,b\right] ;X\right) ;%
\end{array}%
\right.
\end{equation*}%
\begin{equation*}
\leq \left\{ 
\begin{array}{l}
\frac{1}{4}\left( b-a\right) \left\Vert \left\vert f^{\prime }\right\vert
\right\Vert _{\left[ a,b\right] ,1}; \\ 
\\ 
\frac{1}{2^{2+1/q}\left( q+1\right) ^{1/q}}\left( b-a\right)
^{1+1/q}\left\Vert \left\vert f^{\prime }\right\vert \right\Vert _{\left[ a,b%
\right] ,p},\text{ }p>1,\frac{1}{p}+\frac{1}{q}=1,f^{\prime }\in L_{p}\left( %
\left[ a,b\right] ;X\right) ; \\ 
\\ 
\frac{1}{8}\left( b-a\right) ^{2}\left\Vert \left\vert f^{\prime
}\right\vert \right\Vert _{\left[ a,b\right] ,\infty }\text{ }f^{\prime }\in
L_{\infty }\left( \left[ a,b\right] ;X\right) .%
\end{array}%
\right.
\end{equation*}

The inequality (\ref{e4.3}) incorporates \textit{Simpson's rule} as well.
Indeed, if we choose $\alpha =\frac{1}{6},\beta =\frac{4}{6},$ then we get
from (\ref{e4.3}) the result 
\begin{equation}
\left\Vert \left( B\right) \int_{a}^{b}f\left( t\right) dt-\frac{\left(
b-a\right) }{3}\left[ \frac{f\left( a\right) +f\left( b\right) }{2}+2f\left( 
\frac{a+b}{2}\right) \right] \right\Vert  \label{e4.5}
\end{equation}%
\begin{equation*}
\leq \left\{ 
\begin{array}{l}
\frac{1}{3}\left( b-a\right) \left\Vert \left\vert f^{\prime }\right\vert
\right\Vert _{\left[ a,b\right] ,1}; \\ 
\\ 
\frac{1}{\left( q+1\right) ^{1/q}}\frac{\left( 2^{q+1}+1\right) ^{1/q}}{%
6^{1+1/q}}\left( b-a\right) ^{1+1/q}\left\{ \left\Vert \left\vert f^{\prime
}\right\vert \right\Vert _{\left[ a,\frac{a+b}{2}\right] ,p}+\left\Vert
\left\vert f^{\prime }\right\vert \right\Vert _{\left[ \frac{a+b}{2},b\right]
,p}\right\} \\ 
p>1,\frac{1}{p}+\frac{1}{q}=1,f^{\prime }\in L_{p}\left( \left[ a,b\right]
;X\right) ; \\ 
\\ 
\\ 
\frac{5}{72}\left( b-a\right) ^{2}\left[ \left\Vert \left\vert f^{\prime
}\right\vert \right\Vert _{\left[ a,\frac{a+b}{2}\right] ,\infty
}+\left\Vert \left\vert f^{\prime }\right\vert \right\Vert _{\left[ \frac{a+b%
}{2},b\right] ,\infty }\right] ,\text{ }f^{\prime }\in L_{\infty }\left( %
\left[ a,b\right] ;X\right) ;%
\end{array}%
\right.
\end{equation*}%
\begin{equation*}
\leq \left\{ 
\begin{array}{l}
\frac{1}{3}\left( b-a\right) \left\Vert \left\vert f^{\prime }\right\vert
\right\Vert _{\left[ a,b\right] ,1}; \\ 
\\ 
\frac{\left( 2^{q+1}+1\right) ^{1/q}}{2\cdot 3^{1+1/q}\left( q+1\right)
^{1/q}}\left( b-a\right) ^{1+1/q}\left\Vert \left\vert f^{\prime
}\right\vert \right\Vert _{\left[ a,b\right] ,p},\text{ }p>1,\frac{1}{p}+%
\frac{1}{q}=1,f^{\prime }\in L_{p}\left( \left[ a,b\right] ;X\right) ; \\ 
\\ 
\frac{5}{36}\left( b-a\right) ^{2}\left\Vert \left\vert f^{\prime
}\right\vert \right\Vert _{\left[ a,b\right] ,\infty }\text{ }f^{\prime }\in
L_{\infty }\left( \left[ a,b\right] ;X\right) .%
\end{array}%
\right.
\end{equation*}

\begin{remark}
It is obvious that, if the values at $a,\frac{a+b}{2}$ and $b$ of the
function $f:\left[ a,b\right] \rightarrow X$ are available, then one should
choose the rule%
\begin{equation*}
QS\left( f;a,b\right) :=\left( b-a\right) \left[ \frac{f\left( a\right)
+f\left( b\right) }{4}+\frac{1}{2}f\left( \frac{a+b}{2}\right) \right]
\end{equation*}%
that provides a better approximation for the Bochner integral $\left(
B\right) \int_{a}^{b}f\left( t\right) dt$ than the classical Simpson's rule.
\end{remark}

\textbf{2.} If one chooses $x_{1}=\frac{3a+b}{4},x_{2}=\frac{a+b}{2},x_{3}=%
\frac{a+3b}{4},$ then by the use of inequality $\left( \ref{e4.1}\right) ,$
one can derive estimates for the norm of difference%
\begin{equation*}
\left( B\right) \int_{a}^{b}f\left( t\right) dt-\left( b-a\right) \left[
\alpha f\left( \frac{3a+b}{4}\right) +\beta f\left( \frac{a+b}{2}\right)
+\left( 1-\alpha -\beta \right) f\left( \frac{a+3b}{4}\right) \right]
\end{equation*}%
in terms of the Lebesgue norms of the derivative $f^{\prime }.$

We omit the details.

For more scalar-valued three point quadrature rules, see \cite{CD}.

\textbf{Acknowledgement:}{\small \ S. S. Dragomir and Y. J. Cho greatly
acknowledge the financial support from the Brain Pool Program (2002) of the
Korean Federation of Science and Technology Societies. The research was
performed under the "Memorandum of Understanding" between Victoria
University and Gyeongsang National University.}

\end{document}